\documentclass[english]{scrartcl}
\usepackage{babel}

\usepackage{paralist}

\usepackage[lite,author-year]{amsrefs}

\usepackage[utf8]{inputenc}
\usepackage[T1]{fontenc}

\usepackage{amsmath,amssymb}
\DeclareMathOperator{\IM}{Im}
\DeclareMathOperator{\RE}{Re}
\DeclareMathOperator{\co}{co}
\DeclareMathOperator{\dist}{dist}
\DeclareMathOperator{\ran}{ran}
\providecommand{\emb}[1]{i_{#1}}
\providecommand{\field}{\mathbb F}
\providecommand{\naturals}{\mathbb N}
\providecommand{\torus}{\mathbb T}
\providecommand{\bracketcite}[2]{\citeauthor{#1} \ycite{#1}*{#2}}

\usepackage{amsthm}
\newtheorem*{proposition*}{Proposition}
\newtheorem*{theorem*}{Theorem}
\newtheorem{proposition}{Proposition}
\newtheorem{theorem}[proposition]{Theorem}
\theoremstyle{definition}
\newtheorem*{definition*}{Definition}
\theoremstyle{remark}
\newtheorem{remark}{Remark}

\usepackage{mathtools}
\DeclarePairedDelimiter\abs{\lvert}{\rvert}
\DeclarePairedDelimiter\norm{\lVert}{\rVert}
\DeclarePairedDelimiter\set{\lbrace}{\rbrace}
\DeclarePairedDelimiter\dotp{\langle}{\rangle}

\usepackage{hyperref}

\begin{document}
\def\email{pipping@math.fu-berlin.de}
\author{%
  Elias Pipping\\
  \small{\href{mailto:\email}{\nolinkurl{\email}}}
}
\title{L- and M-structure in lush spaces
}
\date{}
\maketitle

{
  \renewcommand\thefootnote{}
}

\begin{abstract}
  Let $X$ be a Banach space which is lush. It is shown that if a
  subspace of $X$ is either an L-summand or an M-ideal then it is also
  lush.
\end{abstract}

\section*{Introduction}
Toeplitz defined \ycite{MR1544315} the \emph{numerical range} of a
square matrix $A$ over the field $\field$ (either $\mathbb R$ or
$\mathbb C$), i.\,e. $A \in \field^{n \times n}$ for some $n \ge 0$,
to be the set
\begin{equation*}
  W(A) = \set{ \dotp{Ax,x} \colon \text{$\norm x = 1$, $x \in \field^n$} },
\end{equation*}
which easily extends to operators on Hilbert spaces. In the 1960s,
Lumer \ycite{MR0133024} and Bauer \ycite{MR0145329} independently
extended this notion to arbitrary Banach spaces. For a Banach space
$X$ whose unit sphere we denote by $S_X$ and an operator $T \in B(X) =
\set{ T \colon X \to X \colon \text{$T$ linear, continuous}}$, we thus
call
\begin{equation*}
  V(T)
  = \set{x^*(Tx) \colon \text{$x^*(x) = 1$, $x^* \in S_{X^*}$, $x \in S_X$}}
  \quad \text{and} \quad
  v(T) = \sup \set{\abs \lambda \colon \lambda \in V(T)}
\end{equation*}
the \emph{numerical range} and \emph{radius} of $T$, respectively. By
construction, we have $v(T) \le \norm T$ for all $T \in B(X)$. The
greatest number $m \ge 0$ that satisfies
\begin{equation*}
  m\norm T \le v(T) \quad \text{for every $T \in B(X)$}
\end{equation*}
is called the \emph{numerical index} of $X$ and denoted by $n(X)$. A
summary of what is and what is not known about the numerical index can
be found in \ocite{MR1823892} and \ocite{MR2267407}. In the special
case $n(X) = 1$ the operator norm and the numerical radius coincide on
$B(X)$.

Several attempts have been made to characterise the spaces with
numerical index one among all Banach spaces geometrically, one of them
in \ocite{MR2296393}. We denote by
\begin{equation*}
  S(B_X, x^*, \alpha) \coloneqq \set{x \in B_X \colon \RE x^*(x) > 1 -\alpha}
\end{equation*}
for any $x^* \in S_{X^*}$ and $\alpha > 0$ an open slice of the unit
ball. Setting $\torus \coloneqq \set{ \omega \in \field \colon \abs
  \omega = 1 } $ and writing $\co(F)$ for the convex hull of a subset
$F \subseteq X$ allows us to write the absolutely convex hull of $F$
as $\co(\torus F)$.
\begin{definition*}
  Let $X$ be a Banach space. If for every two points $u$, $v \in S_X$
  and $\varepsilon > 0$ there is a functional $x^* \in S_{X^*}$ that
  satisfies
  \begin{equation*}
    u \in S(B_X, x^*, \varepsilon)
    \quad \text{and} \quad
    \dist(v, \co(\torus\,S(B_X, x^*, \varepsilon))) < \varepsilon,
  \end{equation*}
  the space $X$ is said to be \emph{lush}.
\end{definition*}
Unfortunately, whilst lush spaces do have numerical index one, spaces
with numerical index one need not be lush \cite{MR2526802}*{Remark
  4.2}. Lushness has proved invaluable in constructing a Banach space
whose dual has strictly smaller numerical index --- answering a
question that up until then had been open for decades. Consequently,
the property deserves attention.

Let us recall some results about sums of Banach spaces.
\begin{proposition*}[\bracketcite{MR1792610}{Proposition~1}]
  Let $(X_n)_{n \in \naturals}$ be a sequence of Banach spaces. Then
  \begin{equation*}
    n\bigl(c_0((X_n)_{n \in \naturals})\bigr)
    = n\bigl(\ell^1((X_n)_{n \in \naturals})\bigr)
    = n\bigl(\ell^\infty((X_n)_{n \in \naturals})\bigr)
    = \inf_{n \in \naturals} n(X_n).
  \end{equation*}
  In particular, the following statements are equivalent:
  \begin{enumerate}[\upshape(i)]
  \item Every $X_n$ has numerical index one,
  \item the space $c_0\bigl((X_n)_{n \in \naturals}\bigr)$ has
    numerical index one,
  \item the space $\ell^1\bigl((X_n)_{n \in \naturals}\bigr)$ has
    numerical index one, and
  \item the space $\ell^\infty\bigl((X_n)_{n \in \naturals}\bigr)$ has
    numerical index one.
  \end{enumerate}
\end{proposition*}

A notion that has been introduced in \ocite{MR0132998} is that of a CL
space. Originally defined for real spaces, it has proven inappropriate
for complex spaces. Thus we will deal with a weakening introduced in
\ocite{MR511813} that had previously been used by
\citeauthor{MR0179580} but remained unnamed.
\begin{definition*}
  Let $X$ be a Banach space. If for every convex subset $F \subseteq
  S_X$ that is maximal in $S_X$ with respect to convexity,
  $\overline{\co}(\torus F) = B_X$ holds, then $X$ is called an
  \emph{almost-CL} space.
\end{definition*}
Almost-CL spaces are easily seen to be lush spaces but the converse
does not hold \cite{MR2296393}*{Example~3.4(c)}. With regard to sums,
the following result has been obtained.
\begin{proposition*}[\bracketcite{MR2056547}{Proposition~8 \& 9}]
  Let $(X_n)_{n \in \naturals}$ be a sequence of Banach spaces. Then
  the following are equivalent:
  \begin{enumerate}[\upshape(i)]
  \item Every $X_n$ is an almost-CL space,
  \item the space $c_0\bigl((X_n)_{n \in \naturals}\bigr)$ is
    almost-CL, and
  \item the space $\ell^1\bigl((X_n)_{n \in \naturals}\bigr)$ is
    almost-CL.
  \end{enumerate}
\end{proposition*}
For the recently introduced lushness property, however, only part of
the corresponding equivalence has been shown.
\begin{proposition*}[\bracketcite{MR2461290}{Proposition~5.3}]
  Let $(X_n)_{n \in \naturals}$ be a sequence of Banach spaces. If
  every $X_n$ is lush, then so are the spaces
  \begin{equation*}
    c_0((X_n)_{n \in \naturals}), \quad%
    \ell^1((X_n)_{n \in \naturals}), \quad \text{and} \quad%
    \ell^\infty((X_n)_{n \in \naturals}).
  \end{equation*}
\end{proposition*}
We seek to improve this result, bringing it up to par with what has
been proved for almost-CL spaces and spaces with numerical index one.
\section*{Inheritance of lushness}
To this end we will show that if $X$ and $Y$ are arbitrary Banach
spaces and one of the two spaces $X \oplus_1 Y$ or $X \oplus_\infty
Y$ is lush, then $X$ and $Y$ are lush themselves.

Such a relation between the spaces $X$, $Y$, and their sum can also be
expressed in terms of projections.
\begin{definition*}
  Let $Z$ be a Banach space and $P \colon Z \to Z$ a linear projection
  that satisfies $\norm z = \max\set{\norm{Pz}, \norm{z - Pz}}$ for
  every $z \in Z$. Then $P$ and $\ran P$ are called an
  \emph{M-projection} and an \emph{M-summand}, respectively.
\end{definition*}
\begin{definition*}
  Let $Z$ be a Banach space and $P \colon Z \to Z$ a linear projection
  that satisfies $\norm z = \norm{Pz} + \norm{z - Pz}$ for every $z
  \in Z$. Then $P$ and $\ran P$ are called an \emph{L-projection} and
  an \emph{L-summand}, respectively.
\end{definition*}
Basic results of L- and M-structure theory that will be used from here
on can be found in \ocite{MR1238713}*{Section~I.1}. If a subspace $X
\subseteq Z$ is an M-summand, its annihilator $X^\bot$ is an L-summand
in $Z^*$. However, an L-summand of $Z^*$ need not be the annihilator
of any space $X \subseteq Z$, nor must subspaces $X \subseteq Z$ for
which $X^\bot$ is an L-summand in $Z^*$ be M-summands. Subspaces $X
\subseteq Z$ for which $X^\bot$ is an L-summand in $Z^*$ are referred
to as \emph{M-ideals}.
\subsection*{M-summands}
We can now proceed to show that M-summands inherit lushness.
\begin{proposition}\label{prop:m-summands-inherit-index-one}
  Let $X$ be an M-summand in a lush space $Z$. Then $X$ is lush.
  \begin{proof}
    Let $u$, $v \in S_X$ and $\varepsilon \in (0,1)$ be arbitrary.
    Since $X$ is an M-summand there is an M-projection $P \colon Z \to
    Z$ with $\ran(P) = X$. Because $Z$ is lush there is a functional
    $z^* \in S_{Z^*}$ satisfying $u \in S(B_Z, z^*, \varepsilon/2)$
    and
    \begin{equation*}
      \dist(v, \co(\torus\,S(B_z, z^*, \varepsilon/2))) < \varepsilon/2.
    \end{equation*}
    Hence there are points $z_1, \dotsc, z_n \in S(B_Z, z^*,
    \varepsilon/2)$ and corresponding $\theta_1, \dotsc, \theta_n \in
    \field$ that satisfy $\sum_{k=1}^n \abs{\theta_k} \le 1$ such that
    $\norm{\sum_{k=1}^n \theta_k z_k - v} < \varepsilon/2$ holds. The
    projection $P$ allows us to split these points up into
    \begin{equation*}
      x_k \coloneqq Pz_k \quad \text{and} \quad y_k \coloneqq P x_k - x_k,
    \end{equation*}
    of which the $x_k$ appear to approximate $v$ mostly by themselves:
    \begin{equation*}
      \norm*{\sum_{k=1}^n \theta_k z_k - v}
      = \max\set*{ \norm*{\sum_{k=1}^n \theta_k y_k},
        \norm*{\sum_{k=1}^n \theta_k x_k - v} }.
    \end{equation*}
    By $\RE z^*(x) > 1 - \varepsilon/2$ and $\norm{z^*} = 1$ we
    clearly have $\RE z^*(y_k) \le \varepsilon/2\norm{x_k} \le
    \varepsilon/2$ for every $k$ and thus
    \begin{equation*}
      \RE z^*(x_k)%
      = \RE z^*(z_k) - \RE z^*(y_k)%
      > 1 - \varepsilon,
    \end{equation*}
    leaving us with $x_k \in S(B_X, z^*, \varepsilon)$, and therefore
    \begin{equation*}
      \dist(v, \co(\torus\,S(B_X,z^*,\varepsilon))) < \varepsilon.
    \end{equation*}
    By restricting $z^*$ to $X$ and normalising the restriction, we
    obtain the desired functional.
  \end{proof}
\end{proposition}
\subsection*{M-ideals}
The celebrated principle of local reflexivity due to \ocite{MR0270119}
can be used to extend
Proposition~\ref{prop:m-summands-inherit-index-one} to M-ideals. More
precisely we require a refined statement.
\begin{theorem*}[\bracketcite{MR0280983}{Section~3}]
  \label{thm:strong-local-reflexivity}
  Let $X$ be a Banach space, $E \subseteq X^{**}$ and $F \subseteq
  X^*$ finite dimensional and $\varepsilon > 0$ arbitrary. Then there
  is an operator $T \colon E \to X$ with $\norm T \norm{T^{-1}} \le 1
  + \varepsilon$ that satisfies $(T \circ \emb{X})(x) = x$ for every
  $x \in X$ with $\emb{X}(x) \in E$ and $x^{**}(x^*) = x^*(Tx^{**})$
  for every $x^* \in F$, $x^{**} \in E$.
\end{theorem*}
An elementary proof is given in \ocite{MR1476378}*{Theorem~2}.
\begin{remark}\label{rem:epsilon-isometries}
  We shall only be concerned with the case $X \ne \set{0}$ in the
  above theorem. Without loss of generality, we can then assume $E
  \cap \emb{X}(X) \ne \set{0}$. Consequently, the
  $\varepsilon$-isometry $T$ can be chosen to satisfy
  \begin{equation*}
    1 - \varepsilon \le \norm{Tz^{**}}
    \le 1 + \varepsilon \quad \text{for every $z^{**} \in S_E$.}
  \end{equation*}
\end{remark}
With that in mind extending
Proposition~\ref{prop:m-summands-inherit-index-one} to M-ideals is
straightforward.
\begin{theorem}
  Let $X$ be an M-ideal in a lush space $Z$. Then $X$ is lush as well.
  \begin{proof}
    Let the points $u$, $v \in S_X$ be arbitrary and $\varepsilon >
    0$. The lushness of $Z$ now guarantees that there is a functional
    $z^* \in S_{Z^*}$ with $u \in S(B_Z, z^*, \varepsilon/2)$ as well
    as an absolutely convex combination of points $z_1, \dotsc, z_n
    \in S(B_Z, z^*, \varepsilon/2)$ and corresponding scalars
    $\theta_1, \dotsc, \theta_n \in \field$ such that
    $\norm{\sum_{k=1}^n \theta_k z_k - v} < \varepsilon/2$ and
    $\sum_{k=1}^n \abs{\theta_k} \le 1$. We observe $Z^{**} =
    X^{\bot\bot} \oplus_\infty M$ for some subspace $M \subseteq
    Z^{**}$. For $k \in \set{1,\dotsc,n}$ we can now find a
    decomposition $\emb{Z}(z_k) = x_k^{**} + y_k^{**}$ with $x_k^{**}
    \in X^{\bot\bot}$ and $y_k^{**} \in M$. By
    \begin{equation*}
      \RE \bigl(\emb{Z^*}(z^*)\bigr)\bigl(\emb{Z}(u)\bigr)%
      = \RE z^*(u)%
      > 1 - \varepsilon/2,
    \end{equation*}
    we clearly have
    \begin{equation*}
      \abs{y^{**}(z^*)}
      \le \varepsilon/2 \quad \text{for every $y^{**} \in S_M$}.
    \end{equation*}
    The functionals $x_k^{**}$ satisfy
    \begin{equation*}
      \RE x_k^{**}(z^*)
      = \RE z^*(z_k) - \RE y_k^{**}(z^*)
      > 1 - \varepsilon
    \end{equation*}
    and in particular
    \begin{equation*}
      1 - \varepsilon \le \norm {x_k^{**}} \le \norm{z_k} = 1.
    \end{equation*}
    We also remark
    \begin{align*}
      \norm*{\sum_{k=1}^n \theta_k z_k - v}%
      &= \max \set*{ \norm*{\sum_{k=1}^n \theta_k y_k^{**}} ,
        \norm*{\sum_{k=1}^n \theta_k x_k^{**} - \emb{Z}(v)} }.
    \end{align*}
    Since $X^{\bot\bot}$ and $X^{**}$ can be identified, we have shown
    that the functionals $x_k^{**}$ meet the requirements of lushness
    for $\emb{X}(u)$ and $\emb{X}(v)$ in $X^{**}$.

    In applying the principle of local reflexivity to the finite
    dimensional subspace $E \coloneqq \operatorname{lin} \set{
      x_1^{**}, \dotsc, x_n^{**}, \emb{Z}(v) } \subseteq X^{**}$, we
    obtain an operator $T \colon E \to X$ that satisfies
    \begin{itemize}
    \item $(T \circ \emb{X})x = x$ for every $x \in X$ with
      $\emb{X}(x) \in E$,
    \item $z^*(Tz^{**}) = z^{**}(z^*)$ for $z^{**} \in E$ and
    \item $1 - \varepsilon/2 \le \norm{Tz^{**}} \le 1 + \varepsilon/2$
      for $z^{**} \in S_E$ (as per
      Remark~\ref{rem:epsilon-isometries}).
    \end{itemize}
    We can now project $x_k^{**}$ onto $X$ with any relevant structure
    preserved. For $x_k \coloneqq Tx_k^{**} \in X$ we observe
    \begin{equation*}
      \norm*{\sum_{k=1}^n \theta_k x_k - v}
      = \norm* {\sum_{k=1}^n \theta_k Tx_k^{**} - (T \circ \emb{Z})v}
      \le (1 + \varepsilon/2) \norm*{\sum_{k=1}^n \theta_k x_k^{**} -
        \emb{Z}(v)} < \varepsilon
    \end{equation*}
    and $\RE z^*(x_k) = \RE x_k^{**}(z^*) > 1 - \varepsilon$. What
    remains to be done is normalising. We thus continue to set $\tilde
    x_k \coloneqq x_k/\norm{x_k}$ and obtain
    \begin{align*}
      \norm{x_k - \tilde x_k}
      &= \abs{\norm{x_k} - 1}\\
      &\le \abs{\norm{x_k} - \norm{x_k^{**}}} + \abs{\norm{x_k^{**}} - 1}\\
      &\le \abs{\norm{Tx_k^{**}} - \norm{x_k^{**}}} + \varepsilon/2\\
      &= \varepsilon \norm{x_k^{**}}/2 + \varepsilon/2\\
      &\le \varepsilon,
    \end{align*}
    and therefore
    \begin{equation*}
      \norm*{\sum_{k=1}^n \theta_k \tilde x_k - v}%
      \le \norm*{\sum_{k=1}^n \theta_k (x_k - \tilde x_k)}
      + \norm*{\sum_{k=1}^n \theta_k x_k - v}
      \le \max_{k \le n} \norm{x_k - \tilde x_k} + \varepsilon
      \le 2 \varepsilon
    \end{equation*}
    as well as
    \begin{equation*}
      \RE z^*(\tilde x_k)
      \ge \RE z^*(x_k) - \norm{x_k - \tilde x_k}
      > 1 - 2 \varepsilon.
    \end{equation*}
  \end{proof}
\end{theorem}
\subsection*{L-summands}
Lushness is also inherited by L-summands. To see this we replace the
complementary parts $y_k$ of $z_k$ with elements $\xi_k \in X$ on
which the functional $z^*$ nearly attains its norm, such that the
$\theta_k \xi_k$ nearly add up to zero.
\begin{theorem}
  Let $X$ be an L-summand of a lush space $Z$. Then $X$ is lush.
  \begin{proof}
    Let $u$, $v \in S_X$ and $\varepsilon > 0$ be arbitrary. Since $Z$
    is lush, for any $\eta > 0$ there is a functional $z^* \in
    S_{Z^*}$ as well as $z_1, \dotsc, z_n \in S(B_Z,z^*,\eta)$ and
    $\theta_1, \dotsc, \theta_n \in \field$ with $\sum_{k=1}^n
    \abs{\theta_k} \le 1$ satisfying $u \in S(B_Z, z^*, \eta)$ and
    $\norm{\sum_{k=1}^n \theta_k z_k - v} < \eta$. Let $P$ be the
    L-projection onto $X$. We set $x_k \coloneqq Pz_k$, $y_k \coloneqq
    z_k - x_k$ and note
    \begin{equation*}
      \norm*{\sum_{k=1}^n \theta_k z_k - v}
      = \norm*{\sum_{k=1}^n \theta_k x_k -v}
      + \norm*{\sum_{k=1}^n \theta_k y_k}.
    \end{equation*}
    In particular, this gives $\norm{\sum_{k=1}^n \theta_k x_k - v} <
    \eta$ and $\norm{\sum_{k=1}^n \theta_k y_k} < \eta$. Replacing
    $y_k$ with $\xi_k \coloneqq \norm{y_k}/\norm u u$ by setting
    $\tilde x_k \coloneqq x_k + \xi_k$ yields $\norm{\tilde x_k} \le
    \norm{z_k} \le 1$ and
    \begin{align*}
      \RE z^*(\tilde x_k)
      &= \RE z^*(z_k - y_k + \xi_k)\\
      &> (1 - \eta) - \norm{y_k} + (1 - \eta) \norm{y_k}\\
      &= 1 - \eta - \eta \norm{y_k}\\
      &\ge 1 - 2 \eta.
    \end{align*}
    We observe
    \begin{align}\label{ineq:real-z*-nearly-attains-norm-on-yk}
      \RE z^*(y_k)%
      = \RE z^*(z_k) - \RE z^*(x_k)
      \ge (1 - \eta) - \norm{x_k}
      \ge \norm{y_k} - \eta,
    \end{align}
    which we will utilise to prove
    \begin{equation}\label{ineq:imag-z*-small-on-yk}
      \bigl(\IM z^*(y_k)\bigr)^2 \le 2 \norm{y_k}\eta.
    \end{equation}
    Since \eqref{ineq:imag-z*-small-on-yk} trivially holds if
    $\norm{y_k} \le \eta$ is satisfied, we shall assume $\norm{y_k} >
    \eta$, leaving us with
    \begin{align*}
      \bigl(\IM z^*(y_k)\bigr)^2%
      &\le (\RE z^*(y_k))^2 + (\IM z^*(y_k))^2 - (\norm{y_k} - \eta)^2\\
      &= \abs{z^*(y_k)}^2 - \norm{y_k}^2 +2 \norm{y_k}\eta - \eta^2\\
      &\le 2 \norm{y_k}\eta - \eta^2\\
      &< 2\norm{y_k}\eta.
    \end{align*}
    We therefore have
    \begin{align*}
      \abs*{\sum_{k=1}^n \theta_k \RE z^*(y_k)}%
      &= \abs*{\sum_{k=1}^n \theta_k z^*(y_k)
        - i \sum_{k=1}^n \theta_k \IM z^*(y_k)}\\
      &\le \norm*{\sum_{k=1}^n \theta_k y_k}
      + \max_{k\le n} \abs*{\IM z^*(y_k)}\\
      &\le \eta + \max_{k\le n} \sqrt{2\norm*{y_k}\eta}\\
      &\le \eta + 2\sqrt\eta.
    \end{align*}
    Applying \eqref{ineq:real-z*-nearly-attains-norm-on-yk} to
    $\delta_k \coloneqq \norm{y_k} - \RE z^*(y_k)$ yields
    $\abs{\delta_k} \le \eta$; we conclude
    \begin{align*}
      \norm*{\sum_{k=1}^n \theta_k \xi_k}
      &\le \abs*{\sum_{k=1}^n \theta_k \RE z^*(y_k)}
      + \abs*{\sum_{k=1}^n \theta_k \delta_k}\\
      &\le 2\eta + 2\sqrt\eta
    \end{align*}
    and thus
    \begin{align*}
      \norm*{\sum_{k=1}^n \theta_k \tilde x_k - v}
      &= \norm*{\sum_{k=1}^n \theta_k \left(x_k + \xi_k\right) - v}\\
      &\le \norm*{\sum_{k=1}^n \theta_k x_k -v}
      + \norm*{\sum_{k=1}^n \theta_k \xi_k}\\
      &\le 3\eta + 2\sqrt\eta.
    \end{align*}
    Going back and choosing $\eta$ such that $3\eta + 2\sqrt\eta <
    \varepsilon$ and $2 \eta < \varepsilon$ are satisfied yields
    \begin{equation*}
      \RE z^*(\tilde x_k) > 1 - \varepsilon
      \quad \text{for every $k \in \set{1,\dotsc,n}$}
    \end{equation*}
    and
    \begin{equation*}
      \dist(v,\co(\torus\,S(B_X, z^*, \varepsilon))) < \varepsilon
    \end{equation*}
    as desired.
  \end{proof}
\end{theorem}

\end{document}